# Meeting an absolute test information target with optimal number of test items via Grand Canonical Monte Carlo simulation

Stefan M. Filipov [a]
Ivan D. Gospodinov [a*]

[a] Department of Informatics, University of Chemical Technology and Metallurgy,
 "Kl. Ohridski"8 Blvd.,
  Sofia 1756, Bulgaria

[*] Corresponding author
  email: idg2@cornell.edu

**Abstract**

This work studies IRT-based Automated Test Assembly (ATA) of multiple test forms (tests) that meet an absolute target information function, i.e. selecting from an item bank only the tests that have information functions that are at a small distance away from the target. The authors introduce the quantities multiplicity of tests and probability of selecting a test with particular number of items $N$ and distance $E$ from the target. A Grand Canonical Monte Carlo test-assembly algorithm is proposed that selects tests according to this probability. The algorithm allows $N$ to vary during the simulation. This work demonstrates that the number of tests that meet the target depends strongly on $N$. The algorithm is capable of finding tests with small values of $E$ and various values of $N$ depending on the need of the test constructor. Most importantly, it can determine the optimal $N$ for which a maximal number of tests with certain specified small $E$ exists.

**Keywords**: density of tests, optimal number of test items, item potential, Grand Canonical Monte Carlo simulation, automated test assembly, target information function

**Introduction**

One of the goals of every test constructor is creating a test that estimates the ability $\theta \in [-3,3]$ of the individual as precisely as possible. A widely used precision measure is the standard error, which in Item Response Theory (IRT) relates to the test information function. According to the three-parameter logistic model of IRT every item $i$ in the test is characterized by the following item information function (Birnbaum, 1968):

$$I_i(\theta) = \left( a_i \frac{p_i(\theta) - c_i}{1 - c_i} \right)^2 \left( \frac{1 - p_i(\theta)}{p_i(\theta)} \right) \tag{1}$$

where $p_i(\theta) = c_i + (1 - c_i)/[1 + \exp(-a_i(\theta - b_i))]$ is the probability of a correct answer to item $i$, and $a_i, b_i, c_i$ are the discrimination, the difficulty, and the guessing probability of item $i$. This holds for dichotomously scored items, i.e. items that are answered either correctly or incorrectly. The test information function $I(\theta)$ is the sum of the item information functions (1) of all the items in the test:

$$I(\theta) = \sum_{i=1}^{N} I_i(\theta) \qquad (2)$$

Typically, in the process of assembling of multiple equivalent tests, the test constructor is presented with requirements about the standard error $SE(\theta)$ from where the absolute target information function $J(\theta)$ is calculated:

$$J(\theta) = 1/SE^2(\theta) \qquad (3)$$

Then, the test constructor needs to select from the item bank a combination of items (a test) that meets the target, i.e. has $I(\theta)$ equal or close to $J(\theta)$. We call such a test a *solution*. When many individuals are tested, in order to maintain test security, typically the test constructor needs to find many solutions (also called parallel forms). This paper will demonstrate that the number of solutions strongly depends on the number of items in the test $N$. Thus, finding $N$ that corresponds to the largest number of solutions may be of great importance.

There exist a variety of methods for construction of test solutions that meet information targets (Lord, 1980, Luecht, 1998, Adema, 1990, Boekkooi-Timminga, 1987, van der Linden & Adema, 1998, Armstrong, Jones, & Kunce, 1998, Veldkamp, 1999, van der Linden, Veldkamp, & Carlson, 2004). None of the existing methods utilize the strong dependence of the number of solutions on $N$. Even though these methods seek multiple solutions, their choice of value for $N$ is not based on knowledge about which value of $N$ corresponds to maximal number of solutions.

This paper proposes a novel Grand Canonical Monte Carlo (GCMC) test-assembly algorithm with varying $N$. The algorithm finds $N$ that corresponds to the largest number of solutions. This $N$ will be called *optimal*. Once the optimal $N$ is found it can be used as a starting test-length constraint for any of the existing ATA methods. The rest of the paper is organized as follows: first the quantities density of tests, item potential, and test selection probability are introduced; then, the test selection probability is used to formulate the GCMC simulation algorithm; then, the algorithm is employed for finding the optimal *N*; and finally, several verification results are presented.

Since this work focuses on the conceptual side of the problem, the simplest possible scenario is considered, i.e. constructing tests that measure a single trait and have no constraints. Incorporation of multiple traits and constraints (such as content balancing, test time, etc.) into the algorithm is possible but it is left for future work.

**Density of tests and test selection probability**

Consider an unordered selection of *N* different items from a bank containing *M* items. Such a selection constitutes a test (a test form) of *N* items. Mathematically it is an *N*-combination of *M* elements. The number of all possible different tests of *N* items is given by the binomial coefficient:

$$\binom{M}{N} = \frac{M!}{N!(M-N)!} \tag{4}$$

The total number of different tests is:

$$\sum_{N=0}^{M} \binom{M}{N} = 2^M \tag{5}$$

Note, that in this analysis there is one test of zero items (0!=1). For a typical item bank the number of item combinations, given by (4), is a huge number for most values of *N*.

Let a test information target function $J(\theta)$ be given that needs to be met with certain precision by the test. This means that the test information function $I(\theta)$, given by (2), needs to be equal or close to $J(\theta)$. As a measure of distance between $J(\theta)$ and $I(\theta)$ for each particular test we employ:

$$E = \int_{-3}^{3} [I(\theta) - J(\theta)]^2 d\theta \tag{6}$$

This is the least squares distance metric widely used in the optimization practice. A different choice of *E* would not, in general, affect the analysis that follows. We seek tests with *E* equal or close to zero. Even though there might be a large number of such tests, they are still a tiny fraction of the huge number of all possible tests. Hence, selecting tests at random with equal probability and waiting for a solution to come up is a hopeless undertaking. This work derives a non-uniform probability distribution for the individual tests to be selected (visited during a simulation) that results in the highest probability of finding tests with pre-specified values of *N* and *E* (e.g. a value of *E* close to zero). Then, a GCMC simulation is proposed as a quick and elegant way of visiting tests (of varying *N*) with this probability.

Let *dE* be a small distance. Let $\Omega(N, E)dE$ be the number of all tests of *N* items and distance between *E* and *E+dE*. Then $\Omega(N, E)$ is the density of tests (the number of tests per unit distance) at *(N,E)*. Note that

$$\int_0^\infty \Omega(N, E)dE = \binom{M}{N} \qquad (7)$$

For the sake of the derivations that follow $\Omega(N, E)$ will be considered continuous and differentiable when necessary.

Consider a *single* particular test of *N* items and distance *E*. Let *P(N,E)* be the probability that this particular test is selected (visited during a simulation). Since we have chosen *P* to depend on *N* and *E* only, then all tests of equal *N* and equal *E* will have equal probability of being selected (that is, equal probability of being visited during a simulation). The exact dependence of *P* on *N* and *E* will be determined shortly. If we multiply $\Omega(N, E)dE$ by *P(N,E)* we will obtain the following probability *p(N,E)dE* that a test selected from the bank has *N* items and distance between *E* and *E+dE*:

$$p(N, E) = \Omega(N, E)P(N, E) \qquad (8)$$

, where *p(N,E)* is the probability density (probability per unit distance *E*) at *(N,E)*. The probability *P* and the probability density *p* satisfy:

$$\sum_{\substack{all \\ tests}} P(N,E) = \sum_{N=0}^{M} \int_{0}^{\infty} p(N,E)dE = 1 \qquad (9)$$

If one wants to find tests with specific, desired number of items $N^*$ and distance $E^*$ then a function $P(N,E)$ must be found such that the probability density $p(N,E)$ has a maximum at $(N^*,E^*)$. Differentiating (8) with respect to $N$ and $E$ gives:

$$\frac{\partial p}{\partial N} = \Omega P \left( \frac{1}{\Omega} \frac{\partial \Omega}{\partial N} + \frac{1}{P} \frac{\partial P}{\partial N} \right) \qquad (10)$$

$$\frac{\partial p}{\partial E} = \Omega P \left( \frac{1}{\Omega} \frac{\partial \Omega}{\partial E} + \frac{1}{P} \frac{\partial P}{\partial E} \right) \qquad (11)$$

Denote the values of the first terms in the brackets at $(N^*, E^*)$ by $\alpha$ and $\beta$:

$$\alpha = \frac{1}{\Omega} \frac{\partial \Omega}{\partial N} \bigg|_{N^*, E^*} \qquad (12)$$

$$\beta = \frac{1}{\Omega} \frac{\partial \Omega}{\partial E} \bigg|_{N^*, E^*} \qquad (13)$$

A sufficient condition that the partial derivatives (10) and (11) become zero at $(N^*,E^*)$ is that $P(N,E)$ satisfies the following differential equations:

$$\frac{1}{P} \frac{\partial P}{\partial N} = -\alpha \qquad (14)$$

$$\frac{1}{P} \frac{\partial P}{\partial E} = -\beta \qquad (15)$$

Integrating (14) and (15) yields:

$$P(N,E) = P(0,0) \exp(-\alpha N - \beta E). \qquad (16)$$

This is the sought non-uniform probability distribution for selecting (visiting) individual tests that results in highest probability density $p$ at the desired, pre-specified $N^*$ and $E^*$. To summarize: if we choose the probability for selecting (visiting) individual tests to be (16), then the product (8) will have a maximum at $(N^*, E^*)$. The constant $P(0,0)$ can, in principle, be obtained from the normalizing condition (9), but it will not be needed. The constants $\alpha$ and $\beta$ can be obtained by calculating $\Omega$ at and around $(N^*, E^*)$ (see eqns. (12) and (13)).

If we introduce the notation

$$\alpha = -\mu/T, \tag{17}$$

$$\beta = 1/T, \tag{18}$$

then the distribution (16) becomes identical to the Gibbs distribution in statistical physics (Kittel & Kroemer, 1980), also known as the Grand Canonical distribution, where $N$ is the number of particles, $E$ is the energy of the system, $T$ is the (fundamental) temperature, $\mu$ is the chemical potential, and $\Omega$ is the density of states. In the Gibbs distribution $P(N,E)$ is the probability of a particular state of a physical system in thermal and diffusion contact with its surroundings at a temperature $T$ and a chemical potential $\mu$. The normalization constant $P(0,0)$ is usually denoted by $1/Z$, $Z$ being the Grand Canonical partition function. From now on, instead of $\alpha$ and $\beta$ in eqn. (16), we will use $-\mu/T$ and $1/T$. We will call $T$ temperature, $\mu$ item potential, $\Omega$ density of tests, and $P$ test selection probability.

In order to investigate the probability density $p(N,E)$ one needs to find $\Omega(N,E)$ first. Finding $\Omega$ for item banks of about 100 or more items becomes a computationally impossible task. As will be seen shortly, the proposed GCMC simulation avoids estimating $\Omega$. Nevertheless, for the sake of verifying the algorithm and demonstrating the features of $p$, we have calculated $\Omega$ for a bank consisting of 50 items for a particular $J(\theta)$ by performing a brute force Monte Carlo counting (see Figure 1). Details will be given in the Results section. Figure 1 shows that tests with moderate $N$ and $E$ are prevalent. For different item banks and targets $J(\theta)$ the function $\Omega(N,E)$ will be different, but it will have a similar shape. In particular, any fixed-$N$ or fixed-$E$ slice of $\Omega(N,E)$ will have a maximum.

**Figure 1**

The density of tests obtained by brute force Monte Carlo counting for a bank of 50 items and
$$J(\theta) = 1 + 3\exp[(-(\theta-1.5)^2)/1.62]/(0.9\sqrt{2\pi})$$

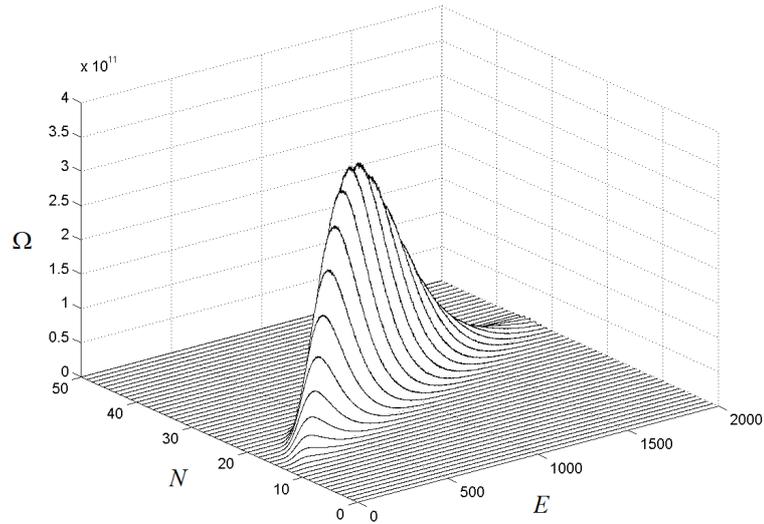

**Figure 2**

Contour plot of $\Omega$, presented in Figure 1, with points (denoted 1 through 9) showing the maxima of nine probability densities $p$ found for different values of $T$ and $\mu$. Details are given in Table 1

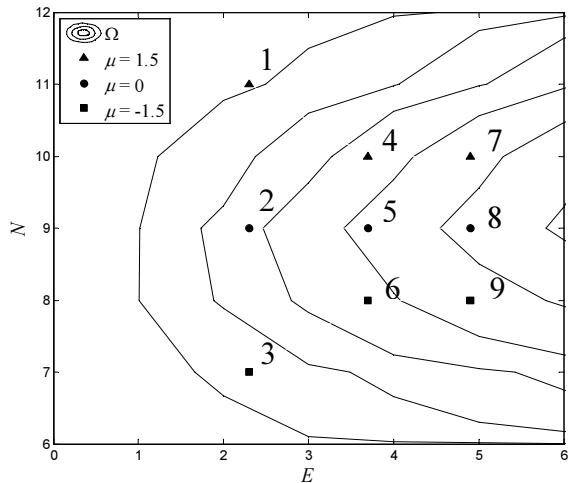

**Table 1**

Data about the nine distributions $p$ with the maxima shown in Figure 2

|  | point 1 | point 2 | point 3 | point 4 | point 5 | point 6 | point 7 | point 8 | point 9 |
|---|---|---|---|---|---|---|---|---|---|
| $T$ | 0.27 | 0.5 | 0.5 | 0.74 | 1 | 1.15 | 1.07 | 1.5 | 1.75 |
| $\Omega$ | 2.27E+02 | 5.15E+04 | 2.95E+03 | 1.41E+05 | 3.68E+05 | 2.18E+05 | 5.44E+05 | 9.91E+05 | 5.07E+05 |
| $s_N$ | 0.8042 | 0.7362 | 0.5500 | 0.9100 | 0.8719 | 0.7900 | 0.9700 | 0.9700 | 0.9200 |
| $s_E$ | 0.9523 | 1.0488 | 1.0100 | 1.9100 | 1.9767 | 2.0500 | 2.5800 | 2.8874 | 3.0000 |

Once $\Omega$ is known, the dependence of $p$ on $T$ and $\mu$ needs to be investigated. Figure 2 shows a contour plot of $\Omega(N,E)$ presented in Figure 1. The region where $E$ approaches zero is of particular interest so only this region is shown in Figure 2. The maxima ($N^*,E^*$) of nine different probability densities $p$, obtained for different values of $T$ and $\mu$, are shown by points denoted 1 through 9. Table 1 shows $T$, $\Omega$, and the standard deviations of $P$ for each point. From eqn. (16), Figures 1 and 2, and Table 1 the following is inferred:

(i)  $p(N,E)$ is sharply peaked around ($N^*,E^*$). It gets narrower as $E^*$ decreases.

(ii) $E^*$ decreases as $T$ decreases. It approaches zero as $T$ approaches zero.

(iii) For any given distance $E^*$ there exists a number of items, called *optimal* and denoted by $N_{opt}$, for which the density of tests $\Omega$ is maximal. If a value $\alpha=0$ (corresponding to $\mu=0$) is substituted into (12) then $[\partial\Omega/\partial N]_{N^*, E^*}=0$ is obtained. Hence, when $\mu=0$ then $N^*=N_{opt}$.

(iv) $N^*$ decreases as $\mu$ decreases. Hence, when $\mu$ is less (respectively more) than zero, then $N^*$ is less (respectively more) than $N_{opt}$.

(v) The density of tests $\Omega$ decreases when $T$ decreases.

Note, that $\Omega$ at point 2, which corresponds to $E^*=2.3$ and $N_{opt}=9$, is 17.5 times greater than $\Omega$ at point 3 and 227 times greater than $\Omega$ at point 1. Results (i) through (v) indicate that the temperature $T$ controls $E^*$ while the item potential $\mu$ controls $N^*$. Result (ii) is not new to the field of Automatic Test Assembly. Several authors (Veldkamp, 1999, van der Linden, Veldkamp, & Carlson, 2004) have used $T$ to control $E^*$. This work introduces the quantities $\mu$ and $\Omega$, derives the probability $p$, defines $N_{opt}$, and proposes a Grand Canonical MC simulation algorithm that, among other things, finds $N_{opt}$.

**Grand Canonical Monte Carlo test-assembly algorithm**

Now that an expression for the individual test probability $P$ is derived (eqn. (16)), an algorithm is needed that ensures that tests are selected from the item bank according to $P$. Thus, most of the time the selected (visited) tests will have values of $N$ and $E$ that are close to $N^*$ and $E^*$ because there lies the maximum of $p$. Since solutions with small $E^*$ are sought, the simulation value for $T$ needs to be small. The choice of a value for $\mu$ depends on the problem that is being solved, e.g. whether one seeks an optimal, a minimal, or a maximal value of $N^*$.

An elegant and efficient way of selecting (visiting) tests with probability (16) is to perform the following GCMC simulation:

1. Specify the values of $T$ and $\mu$.

2. Choose an arbitrary value of $N$ (such that $N \leq M$).

3. Select a combination of $N$ items (a test) from the item bank containing $M$ items. This test is called the *old* test consisting of $N_{old}=N$ items.

4. Calculate $E_{old}$ of this test by using (6).

5. Construct a *new* test from the *old* test by choosing at random one of the following three types of moves:

   - *Replace* a randomly chosen item from the test with a randomly chosen item from the bank. For this move $N_{new}=N_{old}$.

   - *Remove* a randomly chosen item from the test. $N_{new}=N_{old}-1$.

   - *Add* to the test a randomly chosen item from the bank. $N_{new}=N_{old}+1$.

   Removal and addition moves must be chosen with equal probability.

   The item chosen from the bank must be different from those already present in the test.

6. Calculate $E_{new}$ of the new test by using (6).

7. Accept the new test according to the following acceptance probability:

$$acc = \min\left\{1, \binom{M}{N_{new}}\binom{M}{N_{old}}^{-1} \exp\left(\frac{\mu \Delta N - \Delta E}{T}\right)\right\}, \qquad (19)$$

where $\Delta N = N_{new}-N_{old}$, $\Delta E = E_{new}-E_{old}$, and min{} means that one must choose the smaller number from the two given in the curly brackets.

8. If the new test is accepted it becomes the *old* test. If the new test is rejected the old test remains the *old* test. Return to step 5.

Equation (19) is derived using (16) and imposing detailed balance. For particulars about detailed balance the reader is referred to Frenkel & Smit, (2002). Note, that while deriving (19) the quantity $P(0,0)$ gets cancelled. Setting $N_{old}=N$ and using (4) the second term in the curly brackets in (19) becomes:

$$\begin{aligned}
&\exp(-\Delta E/T) && \text{for replacement moves} \\
&(N/(M-N+1))\exp((-\mu-\Delta E)/T) && \text{for removal moves} \quad (20) \\
&((M-N)/(N+1))\exp((\mu-\Delta E)/T) && \text{for addition moves}
\end{aligned}$$

During the GCMC simulation multiple consecutive tests with varying $N$ and $E$ are visited. After some initial equilibration time, $N$ and $E$ start fluctuating around some average values $\langle N \rangle$, $\langle E \rangle$. Since $p$ is sharply peaked, these average values will be very close to $N^*$ and $E^*$.

The algorithm, described in this paper, is principally similar to the GCMC algorithms used in MC molecular simulations (Norman & Filinov, 1969, Frenkel & Smit, 2002), which are proven to be ergodic. Thus, one can be sure that all tests of equal $N$ and equal $E$ appear equally frequently during the simulation. No bias toward any particular solution or group of solutions is introduced.

An important feature of the algorithm is that it can be used for finding the optimal number of items for which there exists a maximum number of solutions. This is achieved by using some small value of $T$ (small enough so that the desired $E^*$ is reached) and a value for the item potential $\mu=0$ (so that $N^*$ is $N_{opt}$). Once $N_{opt}$ is known, the algorithm can be applied for generating multiple solutions with that $N_{opt}$. Alternatively, the value of $N_{opt}$ can be used in another, existing ATA algorithm. Working at $N_{opt}$ may potentially reduce item overlap between the different solutions. Another use of the proposed algorithm is that it can find the minimum or the maximum number of items, depending on the value of $\mu$ used.

The simulation starts from an arbitrary test that is usually far from a solution. If the simulation is started at the desired small value of $T$ the acceptance probability of moves with $\Delta E > 0$ will be very small (see (19)). Thus, often, the system will remain trapped at the starting or some other test with $E$ larger than the desired $E^*$ and will not be able to equilibrate. To avoid this, the method of *simulated annealing* (Kirkpatrick, Gelatt & Vecchi, 1983, van der Linden, Veldkamp & Carlson, 2004, Veldkamp, 1999) is employed where the simulation is started at a higher value of $T$, at which the corresponding distance $E^*$ is higher than the desired, but the acceptance of moves is also higher and the system equilibrates quickly. Then $T$ is lowered gradually, following some cooling schedule, until the desired small value of $E^*$ is reached. At each cooling temperature the system is left to equilibrate and $\langle E \rangle$ is computed to check whether the desired distances $E^*$ has been reached. Finally, the simulation is left to run for some time at the final (lowest) temperature and the average value $\langle N \rangle$ is calculated. As discussed above this average value should be very close to the sought $N^*$. In the case when $\mu=0$ this $N^*$ will be exactly $N_{opt}$.

**Results**

For all the results in this section $E$ was found by evaluating the integrand in eqn. (6) at 61 equidistant points in the interval $\theta \in [-3,3]$. The histogram of the density of tests, shown in Figure 1, was obtained for a bank of 50 items with uniformly distributed parameters $a \in [1,3]$, $b \in [-3,3]$, and constant $c=0.2$. Brute force Monte Carlo counting was used to generate 10 million tests for each $N=1,2,...,50$. The distance interval [0, 2000] was divided into bins of width $dE=1$ and the number of tests that fall into each bin was counted. For each $N$ these numbers were multiplied by the corresponding binomial coefficient and divided by 10 million to obtain $\Omega(N,E)$ (see eqn. (7)).

**Figure 3**
A probability density $p$ for $T=27.41$ and $\mu=7.37$ obtained by: (a) using eqn. (8) and $\Omega$ presented in Figure 1; and (b) performing a GCMC simulation. The means $\langle N \rangle$ and $\langle E \rangle$ are: (a) 15.41 and 130.95; and (b) 15.40 and 131.54. The standard deviations $\sigma_N$ and $\sigma_E$ are: (a) 2.13 and 63.40; and (b) 2.12 and 63.06

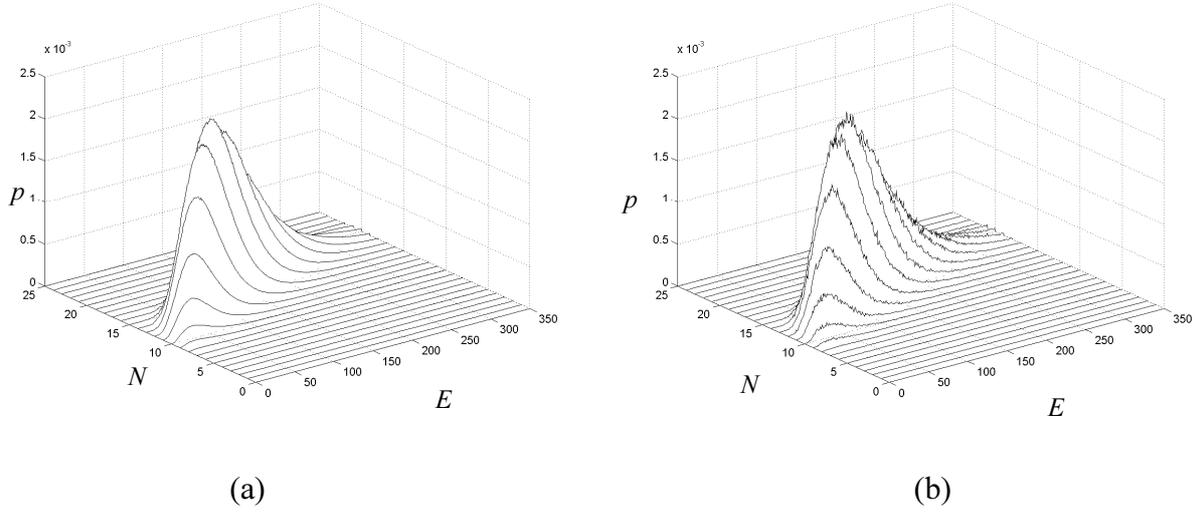

(a)                                                              (b)

The test density $\Omega(N,E)$, shown in Figure 1, was used in order to verify the proposed GCMC algorithm. First, the point $(N^*,E^*)=(15,100)$ was selected in the $(N,E)$ domain. Then, the derivatives of $\Omega$ with respect to $N$ and $E$ at that point were found numerically. Then, using (12) and (13), the values $T=27.41$ and $\mu=7.37$ were found. These values were substituted into (16) to find $P$, then this $P$ was substituted into (8), and finally, after normalizing, the probability density $p$, shown in Figure 3a, was obtained. The same values of $T$ and $\mu$ were used in a GCMC simulation to obtain the probability density $p$ shown in Figure 3b. Both distributions have maxima at $N^*=15$ and $E^*=100$. The mean values and the standard deviations of the two distributions, given in the figure caption, agree very well. Note, that it took 30 hours to calculate $\Omega(N,E)$ needed to obtain the distribution in Figure 3a, while it took less than a minute for the GCMC algorithm to obtain the distribution in Figure 3b.

The GCMC algorithm was applied to a more robust bank of 500 items with uniformly distributed parameters $a \in [1,3]$, $b \in [-3,3]$, and constant $c=0.2$. An absolute target $J(\theta) = 3 + 7\exp[(-(\theta-1.5)^2)/1.62]/(0.9\sqrt{2\pi})$ was imposed. The values of $N^*$ and $\Omega$ corresponding to $E^*=2$ were found for the following five values of $\mu$: -1, -0.5, 0, 0.5, and +1, (see Figure 4). Due to the hugeness of the binomial coefficients the density $\Omega$ could not be determined by brute force Monte Carlo counting. Instead, three simulation distributions $p$ were obtained for $\mu$=-1, 0, and +1, from where the values of $\Omega$ were recovered up to an undetermined factor $C$. The values of $C\Omega$ are given in Figure 4 on a logarithmic scale. Figure 4 shows that the value $N_{opt}$=33, for which $\Omega$ is maximum, is indeed obtained for $\mu$=0. For the other four values of $\mu$, for which $N^*\neq 33$, the density of tests is way smaller, ranging from 25 times to 3 million times smaller. As expected, for a fixed distance $E^*$ the density $\Omega$ has a maximum (sharply peaked) at $\mu$=0 and the number of items $N^*$ decreases as $\mu$ decreases.

**Figure 4**
$N^*$ and $C\Omega$ as functions of $\mu$, obtained by a GCMC simulation using a bank of 500 items for $E^*$=2

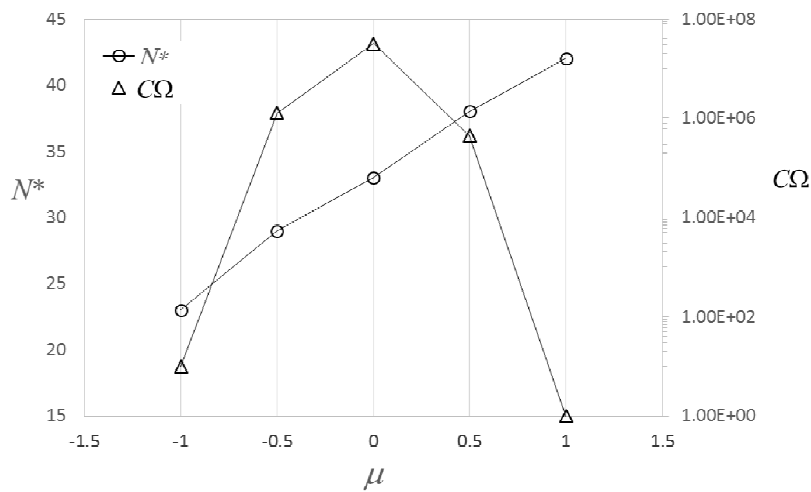

Figure 5a shows the target and the test information functions of one particular solution test with $E^*$=2, $\mu$=0 ($N^*$=33), generated by the GCMC algorithm for the same target and the same item bank that were described in the previous paragraph. Figure 5b shows the histogram of the parameters $a$ and $b$ of all the items in that test. Figures 6a and 6b show analogous distributions of two solution tests with $E^*$=2, $\mu$=-1 ($N^*$=23) and $\mu$=1 ($N^*$=42) (see the corresponding points in Figure 4). Figure 6 demonstrates that tests with less $N^*$ consist of items with higher discriminations $a$. This is to be expected, since items with higher $a$ carry more information than items with lower *a*. For values of $\mu$ less than approximately -1 or more than approximately +1 no solutions with $E^*=2$ exist. This result shows that for every $E^*$ there exists a minimal value of $\mu$ for which $N^*$ is minimum. Of course, this minimum $N^*$ will correspond to a much smaller value of $\Omega$ than the value of $\Omega(N_{opt})$. Finding the exact value of the minimum $N^*$ could be very important in certain test design situations but it is outside the scope of this work.

**Figure 5**

Test and target information functions of a test obtained by GCMC simulation (Figure 5a) with $E^*=2$ and $\mu=0$ ($N^*=33$), and a histogram of the parameters *a* and *b* of all the items in that test (Figure 5b)

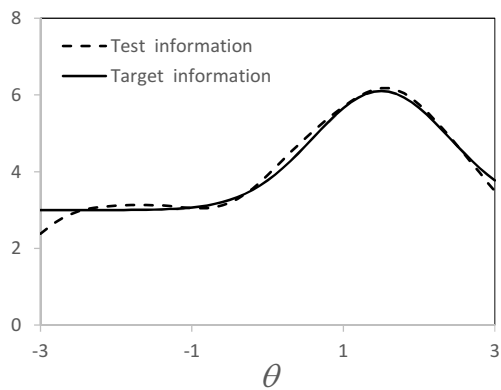 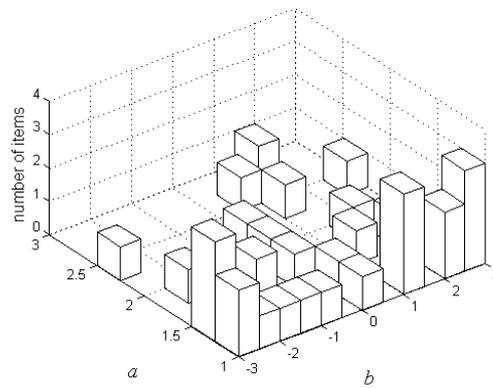

(a) (b)

**Figure 6**

Histogram of the parameters *a* and *b* of all the items in tests obtained by GCMC simulation with $E^*=2$ and (a) $\mu=-1$ ($N^*=23$) and (b) $\mu=1$ ($N^*=42$)

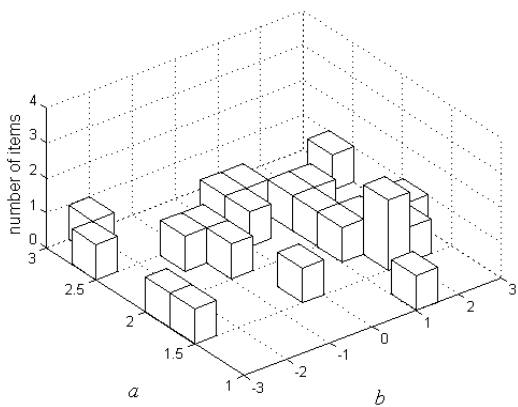 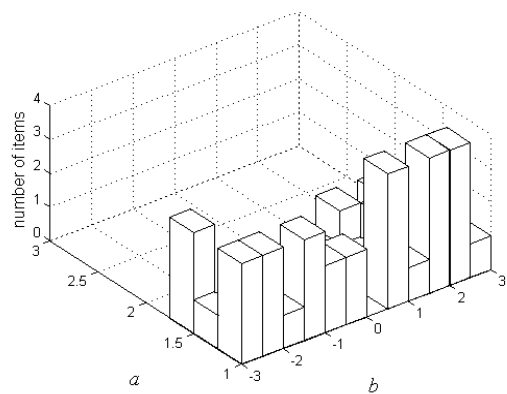

(a) (b)

**Conclusion**

This work introduced the quantities density of tests, item potential, and test selection probability and formulated a simulation algorithm that allows the number of items in the test to vary. The algorithm can be used to solve variety of problems, depending on the values of its control parameters $T$ and $\mu$. The algorithm was used to determine the optimal number of items, that is the number of items that results in highest number of solutions (tests that meet the given target information function). The formalism, introduced in this paper, can serve as a base for future developments of Monte Carlo test-assembly methods, e.g. leading simulations in different ensembles.